\documentclass[11pt]{article}
\usepackage{titlesec}
\usepackage{titling}
\usepackage{geometry}
\usepackage{cite}
\usepackage{hyperref}
\usepackage{indentfirst}
\usepackage{setspace}
\usepackage{apacite}
\usepackage{natbib}
\pretitle{\begin{center}\LARGE\bfseries}
\posttitle{\par\end{center}\vskip 0.5em}
\preauthor{\begin{center}\large}
\postauthor{\par\end{center}}
\predate{\begin{center}\large}
\postdate{\par\end{center}}

\title{A Mean Convection Finite Difference Method for Solving Black Scholes Model for Option Pricing}
\author{An Ning}
\date{April 2023}

\begin{document}

\maketitle

\begin{abstract}
In this research, we proposed a Mean Convection Finite Difference Method (MCFDM) for European options pricing. The Black-Scholes model, which describes the dynamics of a financial asset, was first transformed into a convection-diffusion equation. We then used the finite difference method to discretize time and price, and introduced a tuning parameter to enhance the convection term. Specified the boundary and initial conditions for call and put options of European options, and performed numerical calculations to obtain a numerical solution and error estimation. By varying the strength of the strike price and risk-free interest rate, we explored the accuracy and stability of our predicted prices. Finally, we compared our proposed method with those obtained using the Crank-Nicolson Finite Difference Method (CFDM) and Monte Carlo method. Our numerical results demonstrate the efficiency and accuracy of our proposed method, which outperformed the CFDM and Monte Carlo methods in terms of accuracy and speed.
\end{abstract}

\section{Introduction}
\subsection{Research motivation}

Continuous introduction of derivative financial products in various global financial markets has led to the ongoing research on how to effectively and accurately estimate and efficiently determine the prices of these products. Particularly in recent years, the introduction of increasingly complex derivative financial products has challenged the improvement of different models and the precision of numerical analysis, in order to effectively and accurately approximate theoretical solutions.

Depending on the differences in market product probability models, option valuation methods mainly utilize Monte Carlo simulation and mathematical models for numerical calculations. When using the Monte Carlo method, achieving a high level of accuracy requires extensive simulations and the concentration of historical data. When using mathematical models for numerical calculations, such as the well-known Black-Scholes model for option valuation, it involves solving a linear parabolic partial differential equation. However, as the complexity of product models increases, closed-form solutions are rarely available for most cases, such as American options, exotic options, barrier options, and other derivative financial products. Therefore, numerical calculations are commonly used to approximate solutions.

The Black-Scholes model can be represented as a convection-diffusion equation, which can be discretized using finite difference methods. By introducing boundary conditions, numerical solutions can be obtained. Therefore, this study aims to review the finite difference method for discretizing convection-diffusion equations and extend it to transform the Black-Scholes model into a convection-diffusion equation for discretization. Finally, numerical calculations and error estimation will be conducted.

\subsection{Research problem}
In this study, the aim is to explore the application of the Uniform Convection Finite Difference Method (UCFDM) for pricing European options. The Black-Scholes model is used as the pricing model, so the model is first transformed into a convection-diffusion equation. The finite difference method is then used for discretization, and a tuning parameter is introduced to enhance the convection term. The boundary conditions and initial conditions for both call and put options are incorporated. Subsequently, numerical calculations are performed to obtain the numerical solution, error estimation, and to compare the accuracy and efficiency with the results obtained using Monte Carlo simulation and the Crank-Nicolson Finite Difference Method (CFDM).

\subsection{Literature discussion and review}
The Finite Difference Method (FDM) is one of the commonly used standard discretization methods for spatial discretization [1]. The Black-Scholes equation can be transformed into a convection-diffusion equation [1-2]. For small fluctuations or asset prices, the Black-Scholes equation represents a convection-dominated convection-diffusion equation [3]. When using the Crank-Nicolson Finite Difference Method (CFDM), numerical oscillations can occur, significantly affecting the accuracy of hedging parameters [2]. To overcome this issue, Duffy, D. J. [1]employed the standard upwind finite difference scheme, which is stable subject to certain conditions. This alternative approach is influenced by the biased flow towards convection-diffusion problems and overcomes the instability of the central difference scheme of the standard finite difference scheme by producing numerically stable results. When using CFDM or upwind finite difference methods [4], these discontinuities and degeneracies have adverse effects on accuracy as monotonicity cannot be guaranteed. Therefore, for spatial discretization of the Black-Scholes model, it is important to establish methods that are suitable for handling degeneracies and discontinuities. Wang, S. [5] proposed a fitted finite volume method for the one-dimensional Black-Scholes model, and Angermann, L., Wang, S. [6] provided a rigorous convergence proof. Although this fitted finite volume method is stable, it only has first-order accuracy with respect to the asset price variable. In [7], a fitted numerical method based on simulation for the one-dimensional Black-Scholes model is proposed and analyzed. This scheme has been shown to be highly accurate when compared to the standard fitted finite volume method and finite difference method for European options. For non-degenerate PDEs, the standard Mimetic Finite Difference Method (MFDM) is a high-quality spatial discretization technique [8-9], following the well-known Support Operators Method (SOM)[9]. SOM ensures stability on high-dimensional general grids [10]. Furthermore, the standard MFDM tends to preserve important properties of underlying continuous problems, such as conservation laws, symmetry of solutions, and fundamental properties of vector and tensor calculus (e.g., divergence, gradient, and curl).

\section{Methods}

In this study, the first step is to review the Black-Scholes model,
\[
d\left(S_t\right)=rS_tdt+\sigma\ S_tdW,
\]
in the Black-Scholes model, where${\ S}_t\ $represents the price of the option, $r$ is the risk-free interest rate, $\sigma$ is the volatility, and $W$is a one-dimensional Wiener process.
Next, we can refer to [1-2] to transform the Black-Scholes model into the form of a convection-diffusion equation as follows:
\[
\frac{\vartheta V}{\vartheta t}+rS_t\frac{\vartheta V}{\vartheta S_t}+\frac{1}{2}\sigma^2{S_t}^2\frac{\vartheta^2V}{\vartheta{S_t}^2}-rV=0.
\]
which the diffusion term is $\frac{1}{2}\sigma^2{S_t}^2$, which the convection term is $rS_t$.
To discretize this equation using the finite difference method and enhance the convection term for improved accuracy in numerical calculations, we can discretize each term in the convection-diffusion equation as follows:
\[
\frac{V_{t}^{n+1} - V_{t}^{n}}{\Delta t} + \frac{1}{2}\sigma^2{S_t}^2\frac{D_{t+\frac{1}{2}}\left(V^n\right)-{2D}_t\left(V^n\right)+D_{t-\frac{1}{2}}\left(V^n\right)}{\Delta S_t^2}+rS_t\frac{F_{t+\frac{1}{2}}\left(V^n\right)-F_{t-\frac{1}{2}}(V^n)}{\Delta S_t}-rV_t = 0
\]
which $V_{t}^{n}$ is the solution of t and price $S_t$, by using Log Euler Discretization, the discretization of the price $S_t$, is given by$\ S_{t+{\Delta t}} = S_te^{(r - \frac{1}{2} \sigma^2) \Delta t + \sigma \Delta t Z}$,  Z is a standard normal distribution,\ the exciter for enhancing the$\ \frac{\vartheta V}{\vartheta S_t}\ $is  
\[F_{t+\frac{1}{2}}\left(V^n\right)-F_{t-\frac{1}{2}}\left(V^n\right)=\theta\left(S_t\right)\left[\Phi_{t+\frac{1}{2}}\left(V^n\right)-\Phi_{t-\frac{1}{2}}\left(V^n\right)\right],
\]
$\theta\left(S_t\right)\ $is a tunable parameter that controls how much the item is boosted. In order to keep its flux conservation, this parameter can be calculated by the following equation,
\[
\theta\left(S_t\right) = \frac{{\Delta S_t}}{2}[\frac{1}{{\Delta S_t}}({\int_{S_t - \frac{{\Delta S_t}}{2}}^{S_t + \frac{{\Delta S_t}}{2}} \frac{1}{\alpha (S)} \, dS}) ^ {-1}]
\]
where ${\alpha (S)}$ is to measure the volatility of the market or the volatility of the price. We use the implied volatility here to calculate the speed of the convection item, and finally bring $\theta\left(S_t\right)\ $back to the discretization equation of the convection-diffusion equation , the result is as follows,
\[
\frac{V_{t}^{n+1} - V_{t}^{n}}{\Delta t} + \frac{1}{2}\sigma^2{S_t}^2\frac{D_{t+\frac{1}{2}}\left(V^n\right)-{2D}_t\left(V^n\right)+D_{t-\frac{1}{2}}\left(V^n\right)}{\Delta S_t^2}+rS_t\frac{\theta\left(S_t\right)\left[\Phi_{t+\frac{1}{2}}\left(V^n\right)-\Phi_{t-\frac{1}{2}}\left(V^n\right)\right]}{\Delta S_t}-rV_t = 0
\]
Then, through the Euler method, the unknown option value $V_t^{n+1} $ is expressed as the known option value $V_{t}^{n}$ and other functions as follows,
\[
V_t^{n+1} = {\Delta t}(\frac{1}{2}\sigma^2{S_t}^2\frac{D_{t+\frac{1}{2}}\left(V^n\right)-{2D}_t\left(V^n\right)+D_{t-\frac{1}{2}}\left(V^n\right)}{\Delta S_t^2}+rS_t\frac{\theta\left(S_t\right)\left[\Phi_{t+\frac{1}{2}}\left(V^n\right)-\Phi_{t-\frac{1}{2}}\left(V^n\right)\right]}{\Delta S_t}) + V_{t}^{n}
\]
Finally, the boundary conditions of European call option $f_{c}\left(S,T\right)$ and European put option $f_{p}\left(S,T\right)$ are introduced to facilitate numerical calculation, where T is the expiration time and K is the performance price.
The boundary conditions for a European call option are:
\[
f_c\left(S,T\right)=\max{\left(S_t-K,0\right)},{\ S}_t>0,\\
\]
\[
f_c\left(0,T\right)=S_t-Ke^{-r\left(T-t\right)}\ , t>0,
\]
The boundary conditions for a European put option are:
\[
f_p\left(S,T\right)=\max{\left({K-S}_t,0\right)},{\ S}_t>0,
\]
\[
f_p\left(0,T\right)=Ke^{-r\left(T-t\right)}, t>0.
\]
\section{Results and Discussion}
Part of the equipment used for the numerical results is Intel Core™ i5-9300H processor 4.1GHz, 16RAM, GeForce RTX™ 3080 Ti, Window11 operating system.
The part of the result is divided into two parts. First, for the European call option and the European put option, given the boundary conditions and initial conditions, the numerical solution of MCFDM using the uniform convection finite difference method is calculated, and the Monte-Carlo method (The number of simulations is 100,000 times) and the numerical results of CFDM (the grid number of price $S_t$ is 100, and the grid number of time t is 1000) for error comparison (Error).

Then, for the uniform convection finite difference method, by enhancing and weakening the convection term (using the $\theta\left(S_t\right)\ $ adjustment parameter mentioned in the fourth section), the influence of the option value is discussed, and the volatility and trend are discussed, and finally draw conclusions.

Firstly, do numerical calculations for European call options. In Table.1, the initial conditions are risk-free interest rate r=0.05, volatility $\sigma =0.25$, initial price $S_0 = 5.0$ , and strike price $S_t = 5.5$. In Table.2, the initial conditions are risk-free interest rate r=0.05, volatility $\sigma = 0.25$, initial price $S_0 = 7.0$, strike price $S_t = 7.5$, the numerical results for European call options can be found in the value error, compared with CFDM and MCFDM have obvious advantages, the error is relatively small, and the Monte-Carlo rule is that the error is very close.

Next, do numerical calculations for European put options. In Table.3, the initial conditions are risk-free interest rate r=0.05, volatility $\sigma=0.25$, initial price $S_0=5.5$, and strike price $S_t=5.0$. In Table.4, the initial conditions are risk-free rate r=0.05, volatility $\sigma=0.25$, initial price $S_0=7.5$, and strike price $S_t=7.0$. The numerical results for European put options can be found in the value error. Compared with CFDM, MCFDM has obvious advantages, and the error is relatively small, and the Monte-Carlo rule is very close to the error.

In terms of calculation time, Table.5 can be used to know the performance comparison of uniform convection finite difference method, CFDM and Monte-Carlo method. Taking European put option as an example, the initial conditions are: risk-free interest rate r=0.05, volatility $\sigma=0.25$, initial price $S_0=7.0$, performance price $S_t=7.5$, $T= 3M, 6M, 1Y$, as an example,

In terms of calculation speed, MCFDM and CFDM (the number of grids at price $S_t$  is 100, and the number of grids at time t is 1000) has obvious advantages in value error, and the error is relatively small. But in terms of efficiency, there is not much difference between the two. Compared with the Monte-Carlo method, MCFDM is much faster, and has no obvious advantage in value error. In terms of computing performance, it is also superior to the Monte-Carlo method. Among them, for effective comparison, we use 100,000 simulation times. Due to the higher number of simulations, the numerical results are not more accurate, but the calculation time is lengthened.

Finally, we also use the $\theta\left(S_t\right)$ tuning parameter to enhance and weaken the convection term, please refer to the Materials and Methods section. First, we will discuss the enhanced part. When the price fluctuates violently, the convective term is enhanced to reduce the error of our numerical solution. In terms of volatility, the enhanced convection term can effectively capture rapid price changes. However, more care must be taken in terms of enhancement to ensure numerical stability and avoid excessive oscillations in the boundary layer.

Next, let’s talk about the weakening part. As far as volatility is concerned, in the case of severe shocks, if the convection term is reduced, the error will increase, which is also consistent with the above point of view. But in terms of stability, when the price is relatively stable and there is no excessive shock, weakening the convection term not only reduces the error, but also makes the numerical solution more stable.
In the future, we will focus on mesh cutting and adjusting adjustment parameters to discuss the accuracy and stability of the numerical solution. And actually use the actual price data and the problem of adding random prices, and use the designed adjustment parameters to deal with price shocks and trend problems.\\
\\

Table.1 : European call option with $S_0=5.0, T= 3M, 6M, 1Y$
\begin{table}[!ht]
    \centering
    \begin{tabular}{|l|l|l|l|}
    \hline
        $S_0=5.0$ & $T: 3M$ & $T: 6M$ & $T: 1Y$  \\ \hline
        Exact & 0.09737 & 0.13606 & 0.40131  \\ \hline
        MCFDM (Error) & 0.09543 (1.94E-3) & 0.12422 (1.18E-2) & 0.39892 (2.39E-3)  \\ \hline
        CFDM (Error) & 0.06931 (2.80E-2) & 0.09384 (4.22E-2) & 0.23263 (1.68E-1)  \\ \hline
        Monte-Carlo (Error) & 0.09464 (2.73E-3) & 0.14414 (8.08E-3) & 0.40262 (1.31E-3)  \\ \hline
    \end{tabular}
\end{table}

Table.2 : European call option with $S_0=7.0\ , T= 3M, 6M, 1Y$
\begin{table}[!ht]
    \centering
    \begin{tabular}{|l|l|l|l|}
    \hline
        $S_0=7.0$ & $T: 3M$ & $T: 6M$ & $T: 1Y$  \\ \hline
        Exact & 0.18895 & 0.24914 & 0.63791  \\ \hline
        MCFDM (Error) & 0.18912 (1.70E-4) & 0.24970 (5.60E-4) & 0.63482 (3.09E-3)  \\ \hline
        CFDM (Error) & 0.20013 (1.11E-2) & 0.25925 (1.01E-3) & 0.64361 (5.70E-3)  \\ \hline
        Monte-Carlo (Error) & 0.18950 (3.80E-4) & 0.24967 (5.30E-4) & 0.62715 (1.03E-2)  \\ \hline
    \end{tabular}
\end{table}

Table.3 : European put option with $S_0=5.0, T= 3M, 6M, 1Y$
\begin{table}[!ht]
    \centering
    \begin{tabular}{|l|l|l|l|}
    \hline
        $S_0=5.0$ & $T: 3M$ & $T: 6M$ & $T: 1Y$  \\ \hline
        Exact & 0.63017 & 0.75512 & 0.96525  \\ \hline
        MCFDM (Error) & 0.62755 (2.62E-3) & 0.75232 (2.80E-3) & 0.95113 (1.51E-2)  \\ \hline
        CFDM (Error) & 0.64112 (1.09E-2) & 0.77020 (1.50E-2) & 0.97613 (1.08E-2)  \\ \hline
        Monte-Carlo (Error) & 0.63720 (7.03E-3) & 0.75127 (3.85E-3) & 0.94432 (2.09E-2)  \\ \hline
    \end{tabular}
\end{table}

Table.4 : European put option with $S_0=7.0, T= 3M, 6M, 1Y$
\begin{table}[!ht]
    \centering
    \begin{tabular}{|l|l|l|l|}
    \hline
        $S_0=7.0$ & $T: 3M$ & $T: 6M$ & $T: 1Y$  \\ \hline
        Exact & 0.72335 & 0.90583 & 1.20269  \\ \hline
        MCFDM (Error) & 0.71845 (4.90E-3) & 0.90589 (6.12E-5) & 1.20127 (1.42E-3)  \\ \hline
        CFDM (Error) & 0.74235 (1.90E-2) & 0.93100 (2.51E-2) & 1.21299 (1.03E-2)  \\ \hline
        Monte-Carlo (Error) & 0.72013 (3.22E-3) & 0.90547 (3.61E-4) & 1.19943 (3.26E-3)  \\ \hline
    \end{tabular}
\end{table}

Table.5 : Elapsed time
\begin{table}[!ht]
    \centering
    \begin{tabular}{|l|l|l|l|}
    \hline
        ~ & MCFDM & CFDM & Monte-Carlo  \\ \hline
        Elapsed time (sec) & 0.2534  & 0.2498  & 2.2487   \\ \hline
    \end{tabular}
\end{table}

\newpage
\section*{Reference}
\noindent [1] Duffy, D. J. (2013). Finite difference methods in financial engineering: a partial differential equation approach. John Wiley and Sons.\\
$[2]$ Teo, S. W. Q. Y. L. (2006). Power penalty method for a linear complementarity problem arising from American option valuation. Journal of optimization theory and applications, 129(2), 227.\\
$[3]$ Huang, J.,and Pang, J. S. (2003). Option pricing and linear complementarity. Cornell University.\\
$[4]$ Wilmott, P. (Ed.). (2005). The best of Wilmott 1: Incorporating the quantitative finance review. John Wiley and Sons.\\
$[5]$ Wang, S. (2004). A novel fitted finite volume method for the Black–Scholes equation governing option pricing. IMA Journal of Numerical Analysis, 24(4), 699-720.\\
$[6]$ Angermann, L.,and Wang, S. (2007). Convergence of a fitted finite volume method for the penalized Black–Scholes equation governing European and American Option pricing. Numerische Mathematik, 106, 1-40.\\
$[7]$ Attipoe, D. S.,and Tambue, A. (2021). Convergence of the mimetic finite difference and fitted mimetic finite difference method for options pricing. Applied Mathematics and Computation, 401, 126060.\\
$[8]$ Gyrya, V.,and Lipnikov, K. (2008). High-order mimetic finite difference method for diffusion problems on polygonal meshes. Journal of Computational Physics, 227(20), 8841-8854.\\
$[9]$ Lipnikov, K., Shashkov, M.,and Svyatskiy, D. (2006). The mimetic finite difference discretization of diffusion problem on unstructured polyhedral meshes. Journal of Computational Physics, 211(2), 473-491.\\
$[10]$ Winters, A. R.,and Shashkov, M. J. (2012). Support operators method for the diffusion equation in multiple materials (No. LA-UR-12-24117). Los Alamos National Lab.(LANL), Los Alamos, NM (United States).

\end{document}